\documentclass[10pt]{amsart}
\usepackage{amsfonts}
\usepackage{color}
\usepackage[square,compress,comma, numbers]{natbib}
\usepackage[colorlinks=true, citecolor=blue, linkcolor=blue]{hyperref}
\allowdisplaybreaks[4]
\usepackage{amssymb}
\usepackage{amsmath}
\usepackage{ulem}

\definecolor{c20}{rgb}{0.,0.7,0.}
\definecolor{c30}{rgb}{0.,0.,1.}
\definecolor{c40}{rgb}{1,0.1,0.7}
\definecolor{c50}{rgb}{1,0,0}
\definecolor{c60}{rgb}{1,0.9,0.1}

\allowdisplaybreaks[4]

\newcommand{\abs}[1]{\left\lvert #1 \right\rvert}

\newcommand{\pk}[1]{\mathbb{P} \left \{#1 \right \} }

\newcommand{\BQN}{\begin{eqnarray}}
\newcommand{\EQN}{\end{eqnarray}}
\newcommand{\BQNY}{\begin{eqnarray*}}
\newcommand{\EQNY}{\end{eqnarray*}}

\newcommand{\BS}{\begin{sat}}
\newcommand{\ES}{\end{sat}}
\newcommand{\BT}{\begin{theo}}
\newcommand{\ET}{\end{theo}}
\newcommand{\BL}{\begin{lem}}
\newcommand{\EL}{\end{lem}}
\newcommand{\BK}{\begin{korr}}
\newcommand{\EK}{\end{korr}}

\newcommand{\BD}{\begin{de}}
\newcommand{\ED}{\end{de}}

\newcommand{\BIT}{\begin{itemize}}
\newcommand{\EIT}{\end{itemize}}
\newcommand{\BDI}{\begin{description}}
\newcommand{\EDI}{\end{description}}

\newcommand{\BRM}{\begin{remarks}}
\newcommand{\ERM}{\end{remarks}}

\newcommand{\BEL}{\begin{lem}}
\newcommand{\EEL}{\end{lem}}

\newtheorem{theo}{Theorem}[section]
\newtheorem{sat}[theo]{Proposition}
\newtheorem{de}[theo]{Definition}
\newtheorem{lem}[theo]{Lemma}

\newtheorem{korr}[theo]{Corollary}
\newtheorem{remark}[theo]{Remark}
\newtheorem{remarks}[theo]{Remarks}

\newcommand{\prooftheo}[1]{ \textsc{\bf Proof of Theorem} \ref{#1}:}

\newcommand{\prooflem}[1]{\textsc{\bf Proof of Lemma} \ref{#1}:}

\newcommand{\COM}[1]{}

\newcommand{\QED}{\hfill $\Box$}

\def\rw{\rightarrow}

\def\IF{\infty}

\def\LT{\left}
\def\RT{\right}

\def\rw{\rightarrow}

\def\Var{\text{Var}}

\def\Bu+#1{\mathcal{B}^{\varepsilon+}_{u}(#1)}

\begin{document}

\title[Drawdown and drawup for fractional Brownian motion with trend]{Drawdown and drawup for fractional Brownian motion with trend}

\author{Long Bai}
\address{Long Bai, 
Department of Actuarial Science, 
University of Lausanne\\
UNIL-Dorigny, 1015 Lausanne, Switzerland
}
\email{Long.Bai@unil.ch}



\author{Peng Liu}
\address{Liu, Peng, Department of Actuarial Science, University of Lausanne, UNIL-Dorigny, 1015 Lausanne, Switzerland
}
\email{peng.liu@unil.ch}

\bigskip

\date{\today}
 \maketitle

{\bf Abstract:}
In this paper, we consider the drawdown and drawup of the fractional Brownian motion with trend, 
which corresponds to the logarithm of  geometric fractional Brownian motion representing the stock price in 
financial market. We derive the asymptotics of tail probabilities of the maximum drawdown and maximum drawup  as the 
threshold goes to infinity, respectively. It turns out that the extremes of drawdown leads to new scenarios  of asymptotics
 depending on Hurst index of fractional Brownian motion.

{\bf Key Words:}
Drawdown; Drawup; Fractional Brownian motion; Geometric fractional Brownian motion;  Pickands constant; Piterbarg constant.

{\bf AMS Classification:} Primary 60G15; secondary 60G70

\section{Introduction and Preliminaries}

Drawdown, defined as the distance of present value away from its historical running maximum, is an important indicator of downside risks in financial  risk  management. For instance, the drawdown and the maximum drawdown have
been customarily used as risk measures in finance where they measure the
current drop of a stock price, an index or the value of a portfolio from its running
maximum; see, e.g., \cite{Zhang11, Zhang13}. They  can also be deployed in the context of portfolio optimization as constrains; see, e.g.,\cite{Grossman93, Karatzas95, Chekhlov05, Cherny13, Kardaras17}.  Moreover, drawdown  also arises as reward-to-risk ratio
 in performance measures; see, e.g., \cite{Eling11} for the collections of drawdown-based reward-to-risk ratios. Drawdown processes  also appear 
 in other applications, 
such as applied probability and queueing theory; see, e.g., \cite{Mandjesbook1, Mandjesbook2, Palmowski17, Libin17}. Complementary, drawup,
 the dual of drawdown,  which is  the distance of current value from its historical running minimum, has been encountered in many financial 
 applications; see, e.g., \cite{Pospisil09, Zhang13}.

In the literature, e.g., \cite{RJ03, Guasoni06, SK16}, the stock price  $S$ can be modeled by the
so-called geometric fractional Brownian motion, i.e.,
$$S_t=S_0 \exp\left(\mu t +\sigma B_H(t)-\frac{1}{2}\sigma^2 t^{2H}\right),$$
where $\sigma>0, \mu\in\mathbb{R}$ and $B_H$ is a fractional Brownian motion (fBm) with index $H\in (0,1)$ and covariance function satisfying
$$Cov(B_H(s), B_H(t))=\frac{|s|^{2H}+|t|^{2H}-|s-t|^{2H}}{2}, s,t\geq 0.$$

Note that  $S_t$ is reduced to geometric Brownian motion if $H=1/2$ which has massive applications in Finance.
  To facilitate our analysis, we shall work with the log-prices. This motivates us to  consider the  drawdown and drawup for fBm with trend.
Let $X_t=\sigma B_H(t)-\frac{1}{2}\sigma^2t^{2H}+\mu t,$ $\mu\in \mathbb{R}$. For simplicity, we assume that $\sigma=1.$ The drawdown and drawup processes of  $X$ are defined,  respectively, by
$$D_t=\overline{X}_t-X_t, \quad U_t=X_t-\underline{X}_t,$$
where $\overline{X}_t=\sup_{0\leq s\leq t} X_s$ and $\underline{X}_t=\inf_{0\leq s\leq t}X_s$. For some fixed $T\in(0,\IF)$, we are interested in, for any $u>0$,
\BQN\label{Pro}\pk{\sup_{0\leq t\leq T}D_t>u} \quad \text{and} \quad  \pk{\sup_{0\leq t\leq T}U_t>u}.
\EQN
 Notice that the maximum of drawdown over $[0,T]$ has the interpretation as the largest log-loss up to 
 time $T$ and accordingly, the maximum of drawup can be viewed as the largest log-return; see e.g., \cite{Palmowski17}. Additionally, for $H=\frac{1}{2}$, in context of queueing theory, $D_t$ is the transient queue length process starting at $0$ and the corresponding probability in (\ref{Pro}) represents the  overload probability over $[0,T]$; see, e.g., \cite{Mandjesbook1, Mandjesbook2}.\\
  Note that  for the special case $H=1/2$, the exact expressions of (\ref{Pro})  are obtained in  \cite{Drawdowndis00, Drawdowndis04}; see also \cite{Drawdowndis07} concerning the joint distribution of maximum drawdown and maximum drawup up to an independent exponential time. Due to the fact that fBm is neither a semi-martingale nor a Markov process, the exact expressions for $H\neq \frac{1}{2}$ are not available in literature. Hence in this paper we focus on the asymptotics of (\ref{Pro}) as $u\rw\IF$.\\
It is worthwhile to mention that  infinite series representation  of (\ref{Pro}) in  \cite{Drawdowndis00, Drawdowndis04}   for $H=\frac{1}{2}$  is quite complicated. In contrast, we get concise asymptotics for $H=1/2$ in this paper.
  Theorems \ref{th1} and \ref{th2} in section 2 shows that, for $H=\frac{1}{2}$, as $u\rw\IF$,
   $$\pk{\sup_{0\leq t\leq T}D_t>u}\sim 4\Psi\left(\frac{u+(\mu -\frac{1}{2}T)}{\sqrt{T}}\right), \quad \pk{\sup_{0\leq t\leq T}U_t>u}\sim 4\Psi\left(\frac{u+(\frac{1}{2}T-\mu)}{\sqrt{T}}\right).$$
 The technique used in this paper is  {\it uniform double-sum method}  in \cite{Uniform2017}, which is the development of the so-called {\it double-sum method} widely applied in extreme value theory of Gaussian processes and random fields; see, e.g.,  \cite{Pit96}. As it is shown in Theorem \ref{th1} in section 2, the special trend renders the asymptotics for drawdown  quite different from those of non-centered Gaussian random fields related to fBm in literature (see, e.g., \cite{Pit2001, Pit13, KP2015, KEP20171}), leading to new scenarios of asymptotics according to the value of $H$.\\
 We next introduce some useful notation. We begin with Pickands constant, which is defined by
$$\mathcal{H}_H=\lim_{b\rightarrow\IF}\frac{1}{b}\mathcal{H}_H([0,b])
\quad \text{with }
\mathcal{H}_H ([a,b])=\mathbb{E}\LT\{\sup_{t\in [a,b]}e^{\sqrt{2}B_H(t)-|t|^{2H}}\RT\},  \ \ a<b.
$$
Further, Piterbarg constant is given by, for $\nu>0$,
\BQNY
\mathcal{P}_{H}^{\nu}=\lim_{b\rw\IF} \mathcal{P}_{H}^{\nu}([0,b])\quad \text{with} \quad \mathcal{P}_{H}^{\nu} ([0,b])=\mathbb{E}\LT\{\sup_{t\in [0,b]}e^{\sqrt{2}B_H(t)-(1+\nu)|t|^{2H}}\RT\},\ \ b>0.
\EQNY
We can refer to \cite{Pit96, AdlerTaylor, debicki2002ruin,HP2004,DI2005, Debicki17} for the definition, properties and extensions of Pickands and Piterbarg constants, to \cite{demiro2003simulation, DikerY, Dieker15, GeneralPit16, Harper17} for the bounds and  simulations of Pickands and Piterbarg constants. In particular, by \cite{demiro2003simulation}, we have that
\BQN\label{Piterbarg}
\mathcal{P}_{1/2}^\nu=1+\frac{1}{\nu}, \quad \nu>0.
\EQN
\\
 The organization of paper is as follows. In section 2, the main results are displayed. Section 3 is devoted to the proofs of main theorems in section 2. Proofs of lemmas in section 3 is postponed in Appendix A, followed by some useful lemmas in  Appendix B.

\section{Main Results}\label{s.2}
In this section, we present our main results concerning the asymptotics of (\ref{Pro}) as $u\rw\IF$. In contrast to the infinite series representation in \cite{Drawdowndis00, Drawdowndis04},  the asymptotic expressions in the following theorems are quite concise, which allows us to readily understand the asymptotic behavior of the probability that maximum drawdown ( maximum drawup) exceeds a threshold over finite-time horizon. Let $\Psi(u):=\pk{\mathcal{N}>u}$, with $\mathcal{N}$ a standard normal random variable. Then we have the following results.
\BT\label{th1}
Assume that $0<T<\IF$. \\
If $H>1/2$, then
$$\pk{\sup_{0\leq t\leq T}D_t>u}\sim \Psi\left(\frac{u+\mu T-\frac{1}{2}T^{2H}}{T^H}\right).$$
If $H=1/2$, then
 $$\pk{\sup_{0\leq t\leq T}D_t>u}\sim 4\Psi\left(\frac{u+\mu T-\frac{1}{2}T^{2H}}{T^H}\right).$$
 If $1/4<H<1/2$, then
 $$\pk{\sup_{0\leq t\leq T}D_t>u}\sim \left(H^{-1}2^{-\frac{1}{2H}}T^{2H-1}\mathcal{H}_{H}\right)^2u^{\frac{2}{H}-4}\Psi\left(\frac{u+\mu T-\frac{1}{2}T^{2H}}{T^H}\right).$$
  If $H=1/4$, then
 $$\pk{\sup_{0\leq t\leq T}D_t>u}\sim \left(\mathcal{H}_{\frac{1}{4}}\right)^2T^{-1}\int_0^\IF e^{-x-T^{\frac{1}{4}}x^{\frac{1}{2}}}dx u^{4}\Psi\left(\frac{u+\mu T-\frac{1}{2}T^{2H}}{T^H}\right).$$
  If $0<H<1/4$, then
 $$\pk{\sup_{0\leq t\leq T}D_t>u}\sim H^{-1}2^{-\frac{1}{2H}}T^{2H-2}\Gamma\left(\frac{1}{2H}+1\right)\left(\mathcal{H}_{H}\right)^2 u^{\frac{3}{2H}-2}\Psi\left(\frac{u+\mu T-\frac{1}{2}T^{2H}}{T^H}\right).$$
\ET

\BT\label{th2}
Assume that $0<T<\IF$. \\
If $H>1/2$, then
$$\pk{\sup_{0\leq t\leq T}U_t>u}\sim \Psi\left(\frac{u-\mu T+\frac{1}{2}T^{2H}}{T^H}\right).$$
If $H=1/2$, then
 $$\pk{\sup_{0\leq t\leq T}U_t>u}\sim 4\Psi\left(\frac{u-\mu T+\frac{1}{2}T^{2H}}{T^H}\right).$$
 If $0<H<1/2$, then
 $$\pk{\sup_{0\leq t\leq T}U_t>u}\sim 2^{-\frac{1}{H}-\frac{1}{2}}T^{3H}\sqrt{\frac{\pi}{H^3(H-1)}}\left(\mathcal{H}_{H}\right)^2u^{\frac{2}{H}-3}\Psi\left(\inf_{0\leq s\leq T}\frac{u-\mu (T-s)+\frac{1}{2}(T^{2H}-s^{2H})}{(T-s)^H}\right).$$
\ET

\begin{remark} i) In the extremes of Gaussian processes and random fields associated with fBm for finite-time horizon, e.g.,\cite{Pit2001, Pit13, KP2015, KEP20171}, 
we usually have three different types of asymptotics according to $H$: $H>1/2, H=1/2$ and $H<1/2$. However, Theorem \ref{th1} gives more types of asymptotics due to 
the complexity of the trend that is the combination of linear function ($\mu t$) and power function ($-\frac{1}{2}|t|^{2H}$). As we can see from 
the proof of Theorem \ref{th1}, for $1/4<H<1/2$ only the linear trend contribute to the power part of the asymptotics; for  $H=1/4$, both  linear 
trend and power trend affect  the power part; whereas, for $0<H<1/4$, the power trend has the major influence on the power part of the asymptotics. 
However, this phenomena does not appear in Theorem \ref{th2}, where both of linear trend and power trend contribute to the power part of the asymptotics
 for $0<H<1/2.$\\
ii) We here interpret that the analysis of drawdown and drawup for the case $T=\IF$ is meaningless. Let $T=\IF$ and $\widetilde{B}_H=-B_H$. Then
\BQNY
\sup_{0\leq t<\IF}D_t&=&\sup_{0\leq s\leq t<\IF}\left(B_{H}(s)-B_{H}(t)+\frac{1}{2}(t^{2H}-s^{2H})-\mu(t-s)\right)\\
&=&\sup_{0\leq s\leq t<\IF}\left(\widetilde{B}_{H}(t)-\widetilde{B}_{H}(s)+\frac{1}{2}(t^{2H}-s^{2H})-\mu(t-s)\right)\\
&\geq & \sup_{0\leq s\leq t<\IF}\left(\widetilde{B}_{H}(t)-\widetilde{B}_{H}(s)-(|\mu|+1)(t-s)\right)\\
&=&\sup_{s\geq 0} Q(s),
\EQNY
where $$Q(s)=\sup_{t\geq s}\left(\widetilde{B}_{H}(t)-\widetilde{B}_{H}(s)-(|\mu|+1)(t-s)\right).$$
 Corollary 1 in \cite{Kamil17, Kamil18} shows that  for $H\in (0,1)$
 $$\limsup_{s\rw\IF} \frac{Q(s)}{(\log s)^{\frac{1}{2(1-H)}}}=C>0 \quad a.s..$$
 Therefore we have that for $H\in (0,1)$
 $$\sup_{0\leq t<\IF}D_t\geq \sup_{s\geq 0} Q(s)=\IF \quad a.s..$$
 Note that for $t\geq s\geq 1$ and $H\in(0, 1/2]$,  there exists $C_1>0$ such that
 $$t^{2H}-s^{2H}\leq C_1(t-s).$$
 Hence we can analogously show that for $H\in (0,1/2]$
 \BQNY
 \sup_{0\leq t<\IF}U_t&=&\sup_{0\leq s\leq t<\IF}\left(B_{H}(t)-B_{H}(s)-\frac{1}{2}(t^{2H}-s^{2H})+\mu(t-s)\right)\\
 &\geq&\sup_{1\leq s\leq t<\IF}\left(B_{H}(t)-B_{H}(s)-C_2(t-s)\right)=\IF \quad a.s.,
 \EQNY
 where $C_2$ is a positive constant. We conjecture that for $H>1/2$,
 $$\sup_{0\leq t<\IF}U_t=\IF \quad a.s.$$ also holds,
 which needs more technical analysis similarly to \cite{Kamil17, Kamil18}.
\end{remark}

\section{Proofs}
In this section we give the proof of Theorem \ref{th1}-\ref{th2}. In order to prove the aforementioned theorems, we first present several 
lemmas related to the local behaviors of variance and correlation functions of the underlying Gaussian random fields. In  rest of the paper, denote by 
$Q, Q_i, i=1,2,\dots$ some positive constants that may differ from line to line. Moreover, $$f(u,S,\epsilon)\sim h(u), \quad u\rw\IF, \epsilon\rw 0, S\rw\IF, $$
means that 
$$\lim_{S\rw\IF}\lim_{\epsilon\rw 0}\lim_{u\rw\IF}\frac{f(u,S,\epsilon)}{h(u)}=1.$$
Let
$$\sigma_u^\pm(s,t)=\frac{|t-s|^H}{u\mp\mu(t-s)\pm\frac{1}{2}(t^{2H}-s^{2H})}, \quad 0\leq s\leq t\leq T.$$
$$$$
\BEL\label{lem1}  For $u$ sufficiently large $(0,T)=\arg\sup_{0\leq s \leq t\leq T}\sigma_u^-(s,t)$ is unique and  for any $\delta_u>0$ and $\lim_{u\rw\IF}\delta_u=0$
$$\lim_{u\rw\IF}\sup_{(s,t)\in [0,\delta_u]\times [T-\delta_u, T]}\left|\frac{1-\frac{\sigma_u^-(s,t)}{\sigma_u^-(0,T)}}{\frac{H(T-t)}{T}+\frac{H}{T}s+\frac{1}{2u}s^{2H}}-1\right|=0.$$
\EEL

\BEL\label{lem2} i) For $H\geq \frac{1}{2}$ and $u$ sufficiently large $(0,T)=\arg\sup_{0\leq s \leq t\leq T}\sigma_u^+(s,t)$ is unique and  for any $\delta_u>0$ and $\lim_{u\rw\IF}\delta_u=0$
$$\lim_{u\rw\IF}\sup_{(s,t)\in [0,\delta_u]\times [T-\delta_u, T]}\left|\frac{1-\frac{\sigma_u^+(s,t)}{\sigma_u^+(0,T)}}{\frac{H(T-t)}{T}+\frac{H}{T}s}-1\right|=0.$$
ii) For $0<H<\frac{1}{2}$ and $u$ sufficiently large $(s_u,T)=\arg\sup_{0\leq s \leq t\leq T}\sigma_u^+(s,t)$ is unique and $s_u\sim T^{\frac{1}{1-2H}}u^{-\frac{1}{1-2H}}$. Moreover, for any $\delta_u>0$ and $\lim_{u\rw\IF}\delta_u=0$
$$\lim_{u\rw\IF}\sup_{(s,t)\in [0, s_u+\delta_u]\times [T-\delta_u, T]}\left|\frac{1-\frac{\sigma_u^+(s,t)}{\sigma_u^+(s_u,T)}}{\frac{H(T-t)}{T}+\frac{H(1-H)}{2T^2}(s-s_u)^2}-1\right|=0.$$
\EEL
\BEL\label{lem3} For any $\delta_u>0$ and $\lim_{u\rw\IF}\delta_u=0$
$$\lim_{u\rw\IF}\sup_{(s,t), (s',t')\in [0,\delta_u]\times [T-\delta_u, T]}\left|\frac{1-Corr\left(B_H(t)-B_H(s), B_H(t')-B_H(s')\right)}{\frac{|s-s'|^{2H}+|t-t'|^{2H}}{2T^{2H}}}-1\right|=0.$$
\EEL

\prooftheo{th1}
Observe that
\BQNY
\pk{\sup_{0\leq t\leq T}D_t>u}&=&\pk{\sup_{(s,t)\in A}\left(X_s-X_t\right)>u}\\
&=&\pk{\sup_{(s,t)\in A}\left(B_H(s)-B_H(t)+\mu(s-t)-\frac{1}{2}(s^{2H}-t^{2H})\right)>u}\\
&=&\pk{\sup_{(s,t)\in A}Z_u(s,t)>m(u)},
\EQNY
where $$Z_u(s,t)=\frac{B_H(s)-B_H(t)}{u+\mu(t-s)+\frac{1}{2}(s^{2H}-t^{2H})}m(u), \quad m(u)=\frac{u+\mu T-\frac{1}{2}T^{2H}}{T^H}, \quad A=\{(s,t): 0\leq s\leq t\leq T\}.$$
Thus we have that
\BQN\label{main}
\pk{\sup_{(s,t)\in E_u}Z_u(s,t)>m(u)}&\leq& \pk{\sup_{0\leq t\leq T}D_t>u}\nonumber\\
&\leq& \pk{\sup_{(s,t)\in E_u}Z_u(s,t)>m(u)}+\pk{\sup_{(s,t)\in A\setminus E_u}Z_u(s,t)>m(u)},
\EQN
where $ E_u=[0, (\ln m(u))^2/m^2(u)]\times [T-(\ln m(u))^2/m^2(u), T]$.
In light of Lemma \ref{lem1}, it follows that for $u$ sufficiently large, $\sqrt{Var\left(Z_u(s,t)\right)}=\frac{\sigma_u^-(s,t)}{\sigma_u^-(0,T)}$ attains its maximum over $0\leq s\leq t\leq T$ at unique point $(0,T)$ and
 there exists a positive constant $Q$ such that
$$\sup_{(s,t)\in A\setminus E_u}\sqrt{Var\left(Z_u(s,t)\right)}\leq 1-Q\left(\frac{\ln m(u)}{m(u)}\right)^2.$$
Moreover,
$$\mathbb{E}\left((Z_u(s,t)-Z_u(s',t'))^2\right)\leq Q_1\left(|s-s'|^{2H}+|t-s'|^{2H}\right), \quad (s,t), (s',t')\in A,$$
with $Q_1$ a positive constant.
Hence by Piterbarg Theorem (Theorem 8.1 in \cite{Pit96}), we have for $u$ sufficiently large
\BQN\label{upper}
\pk{\sup_{(s,t)\in A\setminus E_u}Z_u(s,t)>m(u)}\leq Q_2(m(u))^{\frac{2}{H}}\Psi\left(\frac{m(u)}{1-Q\left(\frac{\ln m(u)}{m(u)}\right)^2}\right).
\EQN
Next we analyze $\pk{\sup_{(s,t)\in E_u}Z_u(s,t)>u}$. Let $$\Delta(u)=2^{\frac{1}{2H}}T(m(u))^{-\frac{1}{H}}, \quad E_{u,1}=[0, (\ln m(u))^2/(m^2(u)\Delta(u))]^2.$$ Then rewrite
$$\pk{\sup_{(s,t)\in E_u}Z_u(s,t)>u}=\pk{\sup_{(s,t)\in E_{u,1}}Z_u(\Delta(u)s,T-\Delta(u)t)>u}. $$
We distinguish between $H>\frac{1}{2}$, $H=\frac{1}{2}$, $\frac{1}{4}<H<\frac{1}{2}$, $H=\frac{1}{4}$ and $0<H<\frac{1}{4}$.\\
{\it \underline{Case $H>\frac{1}{2}$}}. In order to apply Lemma \ref{Uniform} in Appendix, we need to check conditions.  By Lemmas \ref{lem1} and \ref{lem3}, we have
\BQN\label{var1}
\lim_{u\rw\IF}\sup_{(s,t)\in E_{u,1}}\left|\frac{1-\sqrt{Var(Z_u(\Delta(u)s,T-\Delta(u)t))}}{\Delta(u)\left(\frac{H}{T}t+\frac{H}{T}s\right)}-1\right|=0,
\EQN
\BQN\label{cor1}
\lim_{u\rw\IF}\sup_{(s,t), (s',t')\in E_{u,1}}\left|m^2(u)\frac{1-Corr\left(Z_u(\Delta(u)s,T-\Delta(u)t), Z_u(\Delta(u)s',T-\Delta(u)t') \right)}{|s-s'|^{2H}+|t-t'|^{2H}}-1\right|=0.
\EQN
These imply that (\ref{vvar3}) and (\ref{ccor3}) hold. Following the notation in Lemma \ref{Uniform}, we have that 
Using the fact that
$$\nu_i=\lim_{u\rw\IF}(m(u))^2\frac{H}{T}\Delta(u)=\IF, i=1,2.$$
Noting that $(0,0)\in E_{u,1}$
and by case iii) in Lemma \ref{Uniform} in Appendix, we have
$$\pk{\sup_{(s,t)\in E_{u,1}}Z_u(\Delta(u)s,T-\Delta(u)t)>m(u)}\sim \Psi(m(u)),$$
which together with (\ref{main}) and (\ref{upper}) establishes the claim.\\
{\it \underline{Case $H=\frac{1}{2}$}}. Note that (\ref{var1}) and (\ref{cor1}) still hold for $H=\frac{1}{2}$. Following the notation in Lemma \ref{Uniform}, we have for $i=1,2,$
$$\nu_i=\lim_{u\rw\IF}(m(u))^2\frac{H}{T}\Delta(u)=2^{\frac{1}{2H}}H=1, \quad \lim_{u\rw\IF} a_i(u)=0, \quad \lim_{u\rw\IF}b_i(u)=\lim_{u\rw\IF} (\ln m(u))^2/(m^2(u)\Delta(u))=\IF.$$
Thus by case ii) in Lemma \ref{Uniform} in Appendix, we have
$$\pk{\sup_{(s,t)\in E_{u,1}}Z_u(\Delta(u)s,T-\Delta(u)t)>m(u)}\sim \left(\mathcal{P}_{1/2}^{1}\right)^2\Psi(m(u)),$$
which combined with (\ref{main}), (\ref{upper})  and (\ref{Piterbarg}) establishes the claim.\\
{\it \underline{Case $\frac{1}{4}<H<\frac{1}{2}$}}.
Let
$$I_{k,l}=[kS, (k+1)S]\times[lS, (l+1)S], k,l\geq 0, \quad N(u)=\left[\frac{(\ln m(u))^2}{m^2(u)\Delta(u)S}\right],$$
$$\Lambda_1(u)=\{(k,l, k', l'): 0\leq k,l,k',l'\leq N(u)+1, I_{k,l}\cap I_{k',l'}\neq \emptyset, (k,l)\neq (k',l')\},$$
$$\Lambda_2(u)=\{(k,l, k', l'): 0\leq k,l,k',l'\leq N(u)+1, I_{k,l}\cap I_{k',l'}= \emptyset\}.$$
Bonferroni inequality gives that
\BQN\label{main3}
\Sigma^-(u)-\Sigma\Sigma_1(u)-\Sigma\Sigma_2(u) \leq \pk{\sup_{(s,t)\in E_{u,1}}Z_u(\Delta(u)s,T-\Delta(u)t)>m(u)}\leq \Sigma^+(u),
\EQN
where
\BQNY
\Sigma^{\pm}(u)&=&\sum _{k,l=0}^{N(u)\pm 1}\pk{\sup_{(s,t)\in I_{k,l}}Z_u(\Delta(u)s,T-\Delta(u)t)>m(u)},\\
\Sigma\Sigma_i(u)&=&\sum_{(k,l,k',l')\in \Lambda_i}\pk{\sup_{(s,t)\in I_{k,l}}Z_u(\Delta(u)s,T-\Delta(u)t)>m(u), \sup_{(s,t)\in I_{k',l'}}Z_u(\Delta(u)s,T-\Delta(u)t)>m(u)}, i=1,2.
\EQNY
{\it \underline{Upper or Lower bounds for $\Sigma^{\pm}(u)$}}.  By Lemma \ref{lem1}, we have
\BQN
\lim_{u\rw\IF}\sup_{(s,t)\in E_{u,1}}\left|\frac{1-\sqrt{Var(Z_u(\Delta(u)s,T-\Delta(u)t))}}{\Delta(u)\frac{H}{T}t+\Delta(u)\frac{H}{T}s+\frac{(\Delta(u))^{2H}}{2u}s^{2H}}-1\right|=0.
\EQN
Thus for any $0<\epsilon<1$, let
$$m_{k,l}^{\pm \epsilon}(u)=m(u)\left(1+(1\pm\epsilon)\left(\Delta(u)\frac{H}{T}\left(l\pm 1\right)S+\Delta(u)\frac{H}{T}\left(k\pm1\right)S+\frac{(\Delta(u))^{2H}}{2u}\left(k\pm 1\right)^{2H}S^{2H}\right)\right).$$
Moreover, denote by
$$Z_{u,k,l}(s,t)=\frac{Z_u(\Delta(u)(kS+s),T-\Delta(u)(lS+t))}{\sqrt{Var\left(Z_u(\Delta(u)(kS+s),T-\Delta(u)(lS+t))\right)}}.$$
Then we have
\BQNY
\Sigma^+(u)&\leq& \sum _{k,l=0}^{N(u)+ 1}\pk{\sup_{(s,t)\in [0,S]^2}Z_{u,k,l}(s,t)>m_{k,l}^{-\epsilon}(u)},\\
\Sigma^-(u)&\geq& \sum _{k,l=0}^{N(u)-1}\pk{\sup_{(s,t)\in [0,S]^2}Z_{u,k,l}(s,t)>m_{k,l}^{+\epsilon}(u)}.
\EQNY
Note that (\ref{cor1}) implies that
\BQN\label{cor2}
\lim_{u\rw\IF}\sup_{(s,t)\in [0,S]^2}\left|(m_{k,l}^{\pm\epsilon}(u))^2\frac{1-Corr\left(Z_{u,k,l}(s,t), Z_{u,k,l}(s',t') \right)}{|s-s'|^{2H}+|t-t'|^{2H}}-1\right|=0.
\EQN
Thus by Lemma \ref{Uniform1}, we have that
\BQNY
\lim_{u\rw\IF}\sup_{0\leq k,l\leq N(u)+1}\left|\frac{\pk{\sup_{(s,t)\in [0,S]^2}Z_{u,k,l}(s,t)>m_{k,l}^{\pm\epsilon}(u)}}{\Psi(m_{k,l}^{\pm\epsilon}(u))}
-(\mathcal{H}_{H}([0,S]))^2\right|=0.
\EQNY
This implies that
\BQN\label{Sigma+}
\Sigma^+(u)&\leq& (\mathcal{H}_{H}([0,S]))^2\sum _{k,l=0}^{N(u)+ 1}\Psi(m_{k,l}^{-\epsilon}(u))\nonumber\\
&\leq& (\mathcal{H}_{H}([0,S]))^2\Psi(m(u))\sum _{k,l=0}^{N(u)+ 1}e^{-(1-\epsilon)\left(m^2(u)\Delta(u)\frac{H}{T}\left(l- 1\right)S+m^2(u)\Delta(u)\frac{H}{T}\left(k-1\right)S+m^2(u)\frac{(\Delta(u))^{2H}}{2u}
\left(k- 1\right)^{2H}S^{2H}\right)}\nonumber\\
&=&: \left(\frac{\mathcal{H}_{H}([0,S])}{S}\right)^2\Psi(m(u))\Theta^-(u,S,\epsilon).
\EQN
and
\BQN\label{Sigma-}
\Sigma^-(u)&\geq& (\mathcal{H}_{H}([0,S]))^2\sum _{k,l=0}^{N(u)- 1}\Psi(m_{k,l}^{+\epsilon}(u))\nonumber\\
&\geq& (\mathcal{H}_{H}([0,S]))^2\Psi(m(u))\sum _{k,l=0}^{N(u)+ 1}e^{-(1+\epsilon)\left(m^2(u)\Delta(u)\frac{H}{T}\left(l+ 1\right)S+m^2(u)\Delta(u)\frac{H}{T}\left(k+1\right)S+m^2(u)\frac{(\Delta(u))^{2H}}{2u}
\left(k+1\right)^{2H}S^{2H}\right)}\nonumber\\
&=&: \left(\frac{\mathcal{H}_{H}([0,S])}{S}\right)^2\Psi(m(u))\Theta^+(u,S,\epsilon).
\EQN
Next we analyze $\Theta^{\pm}(u,S,\epsilon)$.
Note that
$$\sup_{0\leq k\leq N(u)+1}m^2(u)\frac{(\Delta(u))^{2H}}{2u}|k- 1|^{2H}S^{2H}\leq Q (m(u))^{2-4H}\frac{(\ln m(u))^{4H}}{u}\leq Q u^{1-4H}(\ln u)^{4H}\rw  0. $$
Hence, setting 
\BQN\label{upsilon}\upsilon(u, \epsilon)=(1-\epsilon)m^2(u)\Delta(u)\frac{H}{T},
\EQN it follows that
\BQN\label{Theta}
\Theta^-(u,S,\epsilon)&\leq& S^2\sum _{k,l=0}^{N(u)+ 1}e^{-\left(\upsilon(u,\epsilon)\left(l- 1\right)S+\upsilon(u,\epsilon)\left(k-1\right)S\right)}\nonumber\\
&=&\left(\upsilon(u,\epsilon)\right)^{-2}\left(\sum _{l=0}^{N(u)+ 1}e^{-\upsilon(u,\epsilon)\left(l- 1\right)S}\upsilon(u,\epsilon)S\right)
\left(\sum _{k=0}^{N(u)+ 1}e^{-\upsilon(u,\epsilon)\left(k-1\right)S}\upsilon(u,\epsilon)S\right)\nonumber\\
&\leq& \left(\upsilon(u,\epsilon)\right)^{-2}\left(\int_0^\IF e^{-t}dt\right)^2\nonumber\\
&\sim& (m^2(u)\Delta(u))^{-2}\left(\frac{T}{H}\right)^2= \left(H^{-1}2^{-\frac{1}{2H}}T^{2H-1}\right)^2u^{\frac{2}{H}-4}, \quad u\rw\IF, \epsilon\rw 0, S\rw\IF,
\EQN
which together with the fact that
$$\lim_{S\rw\IF}\frac{\mathcal{H}_{H}([0,S])}{S}=\mathcal{H}_H$$
leads to
\BQN\label{upperasym}
\Sigma^+(u)\leq \left(H^{-1}2^{-\frac{1}{2H}}T^{2H-1}\mathcal{H}_H\right)^2u^{\frac{2}{H}-4}\Psi(m(u)), \quad u\rw\IF.
\EQN
Similarly, we can show that
\BQNY
\Theta^+(u,S,\epsilon)\geq \left(H^{-1}2^{-\frac{1}{2H}}T^{2H-1}\right)^2u^{\frac{2}{H}-4}, \quad u\rw\IF, \epsilon\rw 0, S\rw\IF.
\EQNY
Hence
\BQN\label{lowerasym}
\Sigma^-(u)
\geq \left(H^{-1}2^{-\frac{1}{2H}}T^{2H-1}\mathcal{H}_H\right)^2u^{\frac{2}{H}-4}\Psi(m(u)), \quad u\rw\IF.
\EQN
{\it \underline{Upper bounds of $\Sigma\Sigma_i(u), i=1,2$}}. For $(k,l,k',l')\in \Lambda_1$, without loss of generality, we assume that $k'=k+1$. Then denote by
$$I_{k',l'}^{(1)}=[(k+1)S, (k+1)S+\sqrt{S}]\times [l'S,(l'+1)S], \quad I_{k',l'}^{(2)}=[(k+1)S+\sqrt{S}, (k+2)S,]\times [l'S,(l'+1)S].$$
Hence, for $(k,l,k',l')\in \Lambda_1$,
\BQNY
&&\pk{\sup_{(s,t)\in I_{k,l}}Z_u(\Delta(u)s,T-\Delta(u)t)>m(u), \sup_{(s,t)\in I_{k',l'}}Z_u(\Delta(u)s,T-\Delta(u)t)>m(u)}\\
&&\leq \pk{\sup_{(s,t)\in I_{k,l}}\overline{Z}_u(\Delta(u)s,T-\Delta(u)t)>m_{k,l}^{-\epsilon}(u), \sup_{(s,t)\in I_{k',l'}^{(2)}}\overline{Z}_u(\Delta(u)s,T-\Delta(u)t)>m_{k',l'}^{-\epsilon}(u)}\\
&&\quad +\pk{\sup_{(s,t)\in I_{k',l'}^{(1)}}\overline{Z}_{u}(\Delta(u)s,T-\Delta(u)t)>m_{k',l'}^{-\epsilon}(u)},
\EQNY
where
$$\overline{Z}_u(\Delta(u)s,T-\Delta(u)t)=\frac{Z_u(\Delta(u)s,T-\Delta(u)t)}
{\sqrt{Var(Z_u(\Delta(u)s,T-\Delta(u)t))}}.$$

Noting that (\ref{cor2}) holds and
$$\pk{\sup_{(s,t)\in I_{k',l'}^{(1)}}\overline{Z}_{u}(\Delta(u)s,T-\Delta(u)t)>m_{k',l'}^{-\epsilon}(u)}
=\pk{\sup_{(s,t)\in [0,\sqrt{S}]\times[0,S]}Z_{u,k',l'}(s,t)>m_{k',l'}^{-\epsilon}(u)},$$
  by Lemma \ref{Uniform1} in Appendix, we have that
  \BQNY
\lim_{u\rw\IF}\sup_{0\leq k',l'\leq N(u)+1}\left|\frac{\pk{\sup_{(s,t)\in [0,\sqrt{S}]\times [0,S]}Z_{u,k',l'}(s,t)>m_{k,l}^{\pm\epsilon}(u)}}{\Psi(m_{k,l}^{-\epsilon}(u))}
-\mathcal{H}_{H}([0,\sqrt{S}])\mathcal{H}_{H}([0,S])\right|=0.
\EQNY
Using also the fact that $I_{k,l}$ has at most $8$ neighborhoods and
$$\lim_{S\rw\IF}\frac{\mathcal{H}_H([0,\sqrt{S}])}{S}=\lim_{S\rw\IF}
\frac{\mathcal{H}_H([0,\sqrt{S}])}{\sqrt{S}}\lim_{S\rw\IF}S^{-\frac{1}{2}}
=\mathcal{H}_H\lim_{S\rw\IF}S^{-\frac{1}{2}}=0,$$
in light of (\ref{Sigma+}) and (\ref{Theta}), we have
\BQN\label{sum}
&&\sum_{(k,l,k',l')\in \Lambda_1}\pk{\sup_{(s,t)\in I_{k',l'}^{(1)}}\overline{Z}_{u}(\Delta(u)s,T-\Delta(u)t)>m_{k',l'}^{-\epsilon}(u)}\nonumber\\
&&\quad \leq 8\sum_{k',l'=0}^{N(u)+1}\mathcal{H}_{H}([0,\sqrt{S}])\mathcal{H}_{H}([0,S])
\Psi(m_{k,l}^{-\epsilon}(u))\nonumber\\
&&\quad \leq 8\frac{\mathcal{H}([0,\sqrt{S}])}{S}\frac{\mathcal{H}([0,S])}{S}
\left(H^{-1}2^{-\frac{1}{2H}}T^{2H-1}\right)^2u^{\frac{2}{H}-4}\Psi(m(u))\nonumber\\
&&\quad =o\left(u^{\frac{2}{H}-4}\Psi(m(u))\right), \quad u\rw\IF, S\rw \IF.
\EQN
Lemma \ref{lem3} shows that for $u$ sufficiently large and $(s,t), (s',t')\in E_{u,1}$
\BQNY
Corr\left(\overline{Z}_{u}(\Delta(u)s,T-\Delta(u)t), \overline{Z}_{u}(\Delta(u)s',T-\Delta(u)t') \right)>0
\EQNY
and
\BQNY
\lim_{u\rw\IF} \sup_{(s,t)\neq (s',t'),(s,t), (s',t')\in E_{u,1}}\left|(m(u))^2\frac{1-Corr\left(\overline{Z}_{u}(\Delta(u)s,T-\Delta(u)t), \overline{Z}_{u}(\Delta(u)s',T-\Delta(u)t') \right)}{|s-s'|^{2H}+|t-t'|^{2H}}-1\right|=0.
\EQNY
Hence by Lemma \ref{Uniform2} in Appendix, there exists constants $\mathcal{C}, \mathcal{C}_1>0$ such that for $(k,l,k',l')\in \Lambda_1$ and $u$ sufficiently large
\BQNY
&&\pk{\sup_{(s,t)\in I_{k,l}}\overline{Z}_u(\Delta(u)s,T-\Delta(u)t)>m_{k,l}^{-\epsilon}(u), \sup_{(s,t)\in I_{k',l'}^{(2)}}\overline{Z}_u(\Delta(u)s,T-\Delta(u)t)>m_{k',l'}^{-\epsilon}(u)}\\
&&\quad \leq \mathcal{C}S^4 e^{-\mathcal{C}_1 S^{\frac{H}{2}}}\Psi\left(m_{k,l,k',l'}^{-\epsilon}(u)\right);
\EQNY
and for $(k,l,k',l')\in \Lambda_2$ and $u$ sufficiently large
\BQN\label{neg1}
&&\pk{\sup_{(s,t)\in I_{k,l}}\overline{Z}_u(\Delta(u)s,T-\Delta(u)t)>m_{k,l}^{-\epsilon}(u), \sup_{(s,t)\in I_{k',l'}}\overline{Z}_u(\Delta(u)s,T-\Delta(u)t)>m_{k',l'}^{-\epsilon}(u)}\nonumber\\
&&\quad \leq \mathcal{C}S^4 e^{-\mathcal{C}_1 (|k-k'|^2+|l-l'|^2)^{\frac{H}{2}}S^{H}}\Psi\left(m_{k,l,k',l'}^{-\epsilon}(u)\right),
\EQN
where $$m_{k,l,k',l'}^{-\epsilon}(u)=\min(m_{k,l}^{-\epsilon}(u), m_{k',l'}^{-\epsilon}(u)).$$
Consequently, noting that  $I_{k,l}$ has at most $8$ neighborhoods and in light of  (\ref{Sigma+}) and (\ref{Theta})
\BQNY
&&\sum_{(k,l,k',l')\in \Lambda_1}\pk{\sup_{(s,t)\in I_{k,l}}\overline{Z}_u(\Delta(u)s,T-\Delta(u)t)>m_{k,l}^{-\epsilon}(u), \sup_{(s,t)\in I_{k',l'}^{(2)}}\overline{Z}_u(\Delta(u)s,T-\Delta(u)t)>m_{k',l'}^{-\epsilon}(u)}\\
&&\quad \leq  \sum_{(k,l,k',l')\in \Lambda_1} \mathcal{C}S^4 e^{-\mathcal{C}_1 S^{\frac{H}{2}}}\Psi\left(m_{k,l,k',l'}^{-\epsilon}(u)\right)\\
&&\quad \leq  \sum_{(k,l,k',l')\in \Lambda_1} \mathcal{C}S^4 e^{-\mathcal{C}_1 S^{\frac{H}{2}}}\left(\Psi\left(m_{k,l}^{-\epsilon}(u)\right)+
\Psi\left(m_{k',l'}^{-\epsilon}(u)\right)\right)\\
&&\quad \leq \sum_{k,l=0}^{N(u)+1}16\mathcal{C}S^4 e^{-\mathcal{C}_1 S^{\frac{H}{2}}}\Psi\left(m_{k,l}^{-\epsilon}(u)\right) \\
&&\quad \leq QS^2 e^{-\mathcal{C}_1 S^{\frac{H}{2}}} u^{\frac{2}{H}-4}\Psi(m(u))=o\left(u^{\frac{2}{H}-4}\Psi(m(u))\right), \quad u\rw\IF, S\rw \IF.
\EQNY
Therefore, we can conclude that
\BQN\label{ssigma1}
\Sigma\Sigma_1(u)=o\left(u^{\frac{2}{H}-4}\Psi(m(u))\right), \quad u\rw\IF, S\rw \IF.
\EQN
Moreover, by (\ref{neg1}) and (\ref{Sigma+})-(\ref{Theta})
\begin{align}\label{ssigma2}
\Sigma\Sigma_2(u)&\leq \sum_{(k,l,k',l')\in \Lambda_2}\pk{\sup_{(s,t)\in I_{k,l}}\overline{Z}_u(\Delta(u)s,T-\Delta(u)t)>m_{k,l}^{-\epsilon}(u), \sup_{(s,t)\in I_{k',l'}}\overline{Z}_u(\Delta(u)s,T-\Delta(u)t)>m_{k',l'}^{-\epsilon}(u)}\nonumber\\
&\leq \sum_{(k,l,k',l')\in \Lambda_2}\mathcal{C}S^4 e^{-\mathcal{C}_1 (|k-k'|^2+|l-l'|^2)^{\frac{H}{2}}S^{H}}\Psi\left(m_{k,l,k',l'}^{-\epsilon}(u)\right)\nonumber\\
&\leq \sum_{k,l=0}^{N(u)+1}\Psi\left(m_{k,l}^{-\epsilon}(u)\right)2\mathcal{C}S^4\sum_{k',l'\geq 0, k'+l'\neq 0} e^{-\mathcal{C}_1 (|k-k'|^2+|l-l'|^2)^{\frac{H}{2}}S^{H}}\nonumber\\
&\leq \sum_{k,l=0}^{N(u)+1}QS^4e^{-Q_1S^H}\Psi\left(m_{k,l}^{-\epsilon}(u)\right)\nonumber\\
&\leq QS^2e^{-Q_1S^H}u^{\frac{2}{H}-4}\Psi(m(u))=o\left(u^{\frac{2}{H}-4}\Psi(m(u))\right), \quad u\rw\IF, S\rw\IF.
\end{align}
Inserting (\ref{upperasym})-(\ref{lowerasym}) and (\ref{ssigma1})-(\ref{ssigma2}) into (\ref{main3}), we derive that
$$\pk{\sup_{(s,t)\in E_{u,1}}Z_u(\Delta(u)s,T-\Delta(u)t)>m(u)}\sim \left(H^{-1}2^{-\frac{1}{2H}}T^{2H-1}\mathcal{H}_H\right)^2u^{\frac{2}{H}-4}\Psi(m(u)), \quad u\rw\IF,$$
which together with (\ref{main}) and (\ref{upper}) establishes the claim.\\
{\it \underline{Case $H=\frac{1}{4}$}}. Note that (\ref{main3})-(\ref{Sigma-}) still hold for $H=\frac{1}{4}$. We next focus on $\Theta^{\pm}(u,S,\epsilon)$. Recalling that
$$\upsilon(u, \epsilon)=(1-\epsilon)m^2(u)\Delta(u)\frac{H}{T},$$
it follows that
\BQNY
\Theta^-(u,S,\epsilon)&=&S^2\sum _{k,l=0}^{N(u)+ 1}e^{-\left(\upsilon(u, \epsilon)\left(l- 1\right)S+\upsilon(u, \epsilon)\left(k-1\right)S+m^2(u)\frac{(\Delta(u))^{2H}}{2u}\left(k-1\right)^{2H}S^{2H}\right)}\\
&=&\sum _{l=0}^{N(u)+ 1}e^{-\upsilon(u, \epsilon)\left(l-1\right)S} S \sum_{k=0}^{N(u)+1}e^{-\left(\upsilon(u, \epsilon)\left(k-1\right)S+(1-\epsilon)m^2(u)\frac{(\Delta(u))^{2H}}{2u}\left(k-1\right)^{2H}S^{2H}\right)}S.
\EQNY
The first sum satisfies
\BQN\label{firstsum1}
\sum _{l=0}^{N(u)+ 1}e^{-\upsilon(u, \epsilon)\left(l-1\right)S} S
&=&\left(\upsilon(u, \epsilon)\right)^{-1} \sum _{l=0}^{N(u)+ 1}e^{-\upsilon(u, \epsilon)\left(l-1\right)S}\upsilon(u, \epsilon) S\nonumber\\
&\leq& \left(\upsilon(u, \epsilon)\right)^{-1}\int_0^\IF e^{-t}dt\sim \left(m^2(u)\Delta(u)\frac{H}{T}\right)^{-1}, \quad u\rw\IF, \epsilon \rw 0.
\EQN
For the second one
\BQNY
&&\sum_{k=0}^{N(u)+1}e^{-\left(\upsilon(u, \epsilon)\left(k-1\right)S+(1-\epsilon)m^2(u)\frac{(\Delta(u))^{2H}}{2u}\left(k-1\right)^{2H}S^{2H}\right)}S\\
&&\quad=\left(\upsilon(u, \epsilon)\right)^{-1}\sum_{k=0}^{N(u)+1}e^{-\upsilon(u, \epsilon)\left(k-1\right)S+
\left(\frac{(1-\epsilon)^{\frac{1}{2H}}2^{-\frac{1}{2H}}u^{-\frac{1}{2H}}(m(u))^{\frac{1}{H}}\Delta(u)}{\upsilon(u, \epsilon)}\upsilon(u, \epsilon)\left(k-1\right)S\right)^{2H}}\upsilon(u, \epsilon)S.
\EQNY
Note that for $H=\frac{1}{4}$,
$$\frac{(1-\epsilon)^{\frac{1}{2H}}2^{-\frac{1}{2H}}u^{-\frac{1}{2H}}(m(u))^{\frac{1}{H}}\Delta(u)}{\upsilon(u, \epsilon)}\sim \sqrt{T}(1-\epsilon), \quad u\rw\IF.$$
Thus
\BQNY
&&\sum_{k=0}^{N(u)+1}e^{-\left(\upsilon(u, \epsilon)\left(k-1\right)S+(1-\epsilon)m^2(u)\frac{(\Delta(u))^{2H}}{2u}\left(k-1\right)^{2H}S^{2H}\right)}S\\
&&\sim \left(m^2(u)\Delta(u)\frac{H}{T}\right)^{-1}\int_0^\IF e^{-x-T^{\frac{1}{4}}\sqrt{x}}dx, \quad u\rw\IF, \epsilon\rw 0.
\EQNY
Consequently,
\BQNY
\Theta^-(u,S,\epsilon)\leq \left(m^2(u)\Delta(u)\frac{H}{T}\right)^{-2}\int_0^\IF e^{-x-T^{\frac{1}{4}}\sqrt{x}}dx, \quad u\rw\IF, \epsilon\rw 0.
\EQNY
Similarly,
\BQNY
\Theta^+(u,S,\epsilon)\geq \left(m^2(u)\Delta(u)\frac{H}{T}\right)^{-2}\int_0^\IF e^{-x-T^{\frac{1}{4}}\sqrt{x}}dx, \quad u\rw\IF, \epsilon\rw 0.
\EQNY
In light of (\ref{Sigma+}) and (\ref{Sigma-}), we have that
\BQNY
\Sigma^-(u)&\leq& \left(\frac{\mathcal{H}_{H}([0,S])}{S}\right)^2\left(m^2(u)\Delta(u)\frac{H}{T}\right)^{-2}\int_0^\IF e^{-x-T^{\frac{1}{4}}\sqrt{x}}dx\Psi(m(u))\\
&\leq&\left(H^{-1}2^{-\frac{1}{2H}}T^{2H-1}\mathcal{H}_{H}\right)^2\int_0^\IF e^{-x-T^{\frac{1}{4}}\sqrt{x}}dx u^{\frac{2}{H}-4}\Psi(m(u)), \quad u\rw\IF, S\rw\IF,\\
\Sigma^+(u)&\geq&\left(H^{-1}2^{-\frac{1}{2H}}T^{2H-1}\mathcal{H}_{H}\right)^2\int_0^\IF e^{-x-T^{\frac{1}{4}}\sqrt{x}}dx u^{\frac{2}{H}-4}\Psi(m(u)), \quad u\rw\IF, S\rw\IF.
\EQNY
The negligibility of $\Sigma\Sigma_i(u), i=1,2$ holds due to the fact that (\ref{sum})-(\ref{ssigma2}) are also valid  for $H=\frac{1}{4}$. Therefore we have
$$\pk{\sup_{(s,t)\in E_{u,1}}Z_u(\Delta(u)s,T-\Delta(u)t)>m(u)}\sim \left(H^{-1}2^{-\frac{1}{2H}}T^{2H-1}\mathcal{H}_{H}\right)^2\int_0^\IF e^{-x-T^{\frac{1}{4}}\sqrt{x}}dx u^{\frac{2}{H}-4}\Psi(m(u)), \quad u\rw\IF,$$
which combined with (\ref{main}) and (\ref{upper}) establishes the claim.\\
{\it \underline{Case $0<H<\frac{1}{4}$}}.
For $0<H<\frac{1}{4}$, (\ref{main3})-(\ref{Sigma-}) are satisfied.  In order to get the upper or lower bounds of $\Sigma^{\pm}(u)$, it suffices to analyze $\Theta^{\pm}(u,S,\epsilon)$. Denote by
$$\upsilon'(u, \epsilon)=(1-\epsilon)^{\frac{1}{2H}}2^{-\frac{1}{2H}}
u^{-\frac{1}{2H}}(m(u))^{\frac{1}{H}}\Delta(u),$$
it follows that
\BQNY
\Theta^-(u,S,\epsilon)&=&S^2\sum _{k,l=0}^{N(u)+ 1}e^{-(1-\epsilon)\left(m^2(u)\Delta(u)\frac{H}{T}\left(l- 1\right)S+m^2(u)\Delta(u)\frac{H}{T}\left(k-1\right)S+m^2(u)\frac{(\Delta(u))^{2H}}{2u}\left(k-1\right)^{2H}S^{2H}\right)}\\
&=&\sum _{l=0}^{N(u)+ 1}e^{-\upsilon(u,\epsilon)\left(l-1\right)S} S \sum_{k=0}^{N(u)+1}e^{-\left(\upsilon(u,\epsilon)\left(k-1\right)S+\left(\upsilon'(u,\epsilon)(k-1)S\right)^{2H}\right)}S,
\EQNY
where $\upsilon(u,\epsilon)$ is defined in (\ref{upsilon}).
The first sum satisfies (\ref{firstsum1}) with $0<H<1/4$.
For the second sum
\BQNY
&&\sum_{k=0}^{N(u)+1}e^{-(\upsilon(u,\epsilon)\left(k-1\right)S+\left(\upsilon'(u,\epsilon)(k-1)S\right)^{2H})}S\\
&&=\left(\upsilon'(u,\epsilon)\right)^{-1}
 \sum_{k=0}^{N(u)+1}
e^{-\gamma(u)\upsilon'(u,\epsilon)\left(k-1\right)S+
\left(\upsilon'(u,\epsilon)\left(k-1\right)S\right)^{2H}}
\upsilon'(u,\epsilon)S,
\EQNY
where
\BQNY
\gamma(u)&=&\frac{\upsilon(u,\epsilon)}{\upsilon'(u,\epsilon)}\\
&=&\frac{(1-\epsilon)m^2(u)\Delta(u)\frac{H}{T}}{(1-\epsilon)^{\frac{1}{2H}}2^{-\frac{1}{2H}}u^{-\frac{1}{2H}}(m(u))^{\frac{1}{H}}\Delta(u)}\\
&\sim& Qu^{2-\frac{1}{2H}}\rw 0, \quad u\rw\IF.
\EQNY
Thus
\BQNY
&&\sum_{k=0}^{N(u)+1}e^{-(\upsilon(u,\epsilon)\left(k-1\right)S+\left(\upsilon'(u,\epsilon)(k-1)S\right)^{2H})}S\\
&&\quad \sim \left(\upsilon'(u,\epsilon)\right)^{-1}\int_0^\IF e^{-x^{2H}}dx\\
&&\quad \sim \Gamma\left(\frac{1}{2H}+1\right)T^{-1}u^{\frac{1}{2H}}, \quad u\rw\IF, \epsilon\rw 0.
\EQNY
Consequently,
\BQNY
\Theta^-(u,S,\epsilon)\leq H^{-1}2^{-\frac{1}{2H}}T^{2H-2}\Gamma\left(\frac{1}{2H}+1\right) u^{\frac{3}{2H}-2}, \quad u\rw\IF, \epsilon\rw 0.
\EQNY
Similarly,
\BQNY
\Theta^+(u,S,\epsilon)\geq H^{-1}2^{-\frac{1}{2H}}T^{2H-2}\Gamma\left(\frac{1}{2H}+1\right) u^{\frac{3}{2H}-2}, \quad u\rw\IF, \epsilon\rw 0.
\EQNY
In light of (\ref{Sigma+}) and (\ref{Sigma-}), we have that, as $u\rw\IF, S\rw\IF$,
\BQNY
\Sigma^-(u)&\leq& H^{-1}2^{-\frac{1}{2H}}T^{2H-2}\Gamma\left(\frac{1}{2H}+1\right)\left(\mathcal{H}_{H}\right)^2 u^{\frac{3}{2H}-2}\Psi(m(u)),\\
\Sigma^+(u)&\geq&H^{-1}2^{-\frac{1}{2H}}T^{2H-2}\Gamma\left(\frac{1}{2H}+1\right)\left(\mathcal{H}_{H}\right)^2 u^{\frac{3}{2H}-2}\Psi(m(u)).
\EQNY
Following line by line the same as (\ref{sum})-(\ref{ssigma2}), we can show that for $i=1,2$
$$\Sigma\Sigma_i(u)=o\left(u^{\frac{3}{2H}-2}\Psi(m(u))\right), \quad u\rw\IF, S\rw\IF.$$
Therefore, we conclude that
$$\pk{\sup_{(s,t)\in E_{u,1}}Z_u(\Delta(u)s,T-\Delta(u)t)>m(u)}\sim H^{-1}2^{-\frac{1}{2H}}T^{2H-2}\Gamma\left(\frac{1}{2H}+1\right)\left(\mathcal{H}_{H}\right)^2 u^{\frac{3}{2H}-2}\Psi(m(u)), \quad u\rw\IF,$$
which establishes the claim with aid of  (\ref{main}) and (\ref{upper}). This completes the proof. \QED

\prooftheo{th2} We distinguish between $H\geq \frac{1}{2}$ and $H<\frac{1}{2}$.\\
{\it \underline{Case $H\geq \frac{1}{2}$}}. We have that
\BQNY
\pk{\sup_{0\leq t\leq T}U_t>u}&=&\pk{\sup_{(s,t)\in A} (X_t-X_s)>u}\\
&=&\pk{\sup_{(s,t)\in A}(B_H(t)-B_H(s)-\frac{1}{2}(t^{2H}-s^{2H})+\mu(t-s))>u}\\
&=&\pk{\sup_{(s,t)\in A} Z_{u,1}(s,t)>m_1(u)},
\EQNY
where $$Z_{u,1}(s,t)=\frac{B_H(t)-B_H(s)}{u-\mu(t-s)+\frac{1}{2}(t^{2H}-s^{2H})}m_1(u), \quad m_1(u)=\frac{u-\mu T+\frac{1}{2}T^{2H}}{T^H}, \quad A=\{(s,t), 0\leq s\leq t\leq T\}.$$
Furthermore,
\BQN\label{main1}
\pk{\sup_{(s,t)\in E_{u,2}} Z_{u,1}(s,t)>m_1(u)}&\leq& \pk{\sup_{0\leq t\leq T}U_t>u}\nonumber\\
&\leq& \pk{\sup_{(s,t)\in E_{u,2}} Z_{u,1}(s,t)>m_1(u)}+\pk{\sup_{(s,t)\in A\setminus E_{u,2}} Z_{u,1}(s,t)>m_1(u)},
\EQN
where $$E_{u,2}=[0, (\ln m_1(u))^2/(m_1(u))^2]\times[T-(\ln m_1(u))^2/(m_1(u))^2,T].$$
In light of Lemma \ref{lem2}, it follows that for $u$ sufficiently large
\BQNY
\sup_{(s,t)\in A\setminus E_{u,2}}\sqrt{Var(Z_{u,1}(s,t))}\leq 1-Q\left(\frac{\ln m_1(u)}{m_1(u)}\right)^2.
\EQNY
Moreover, direct calculation shows that
\BQNY
\mathbb{E}\left((Z_{u,1}(s,t)-Z_{u,1}(s',t'))^2\right)\leq Q_1(|t-t'|^{2H}+|s-s'|^{2H}), \quad (s,t), (s',t')\in A.
\EQNY
 Using Piterbarg Theorem (Theorem 8.1 in \cite{Pit96}), we have for $u$ sufficiently large
 \BQN\label{upper1}
 \pk{\sup_{(s,t)\in A\setminus E_{u,2}} Z_{u,1}(s,t)>m_1(u)}\leq Q_2 (m_1(u))^{\frac{2}{H}}\Psi\left(\frac{m_1(u)}{1-Q\left(\frac{\ln m_1(u)}{m_1(u)}\right)^2}\right).
 \EQN
 Next we focus on $\pk{\sup_{(s,t)\in E_{u,2}} Z_{u,1}(s,t)>m_1(u)}$.  Lemmas \ref{lem1} and \ref{lem3} lead to
 \BQNY
\lim_{u\rw\IF}\sup_{(s,t)\in E_{u,2}}\left|\frac{1-\sqrt{Var(Z_{u,1}(s,t))}}{\frac{H(T-t)}{T}+\frac{H}{T}s}-1\right|=0,\quad \lim_{u\rw\IF}\sup_{(s,t), (s',t')\in E_{u,2}}\left|\frac{1-Corr\left(Z_{u,1}(s,t), Z_{u,1}(s',t')\right)}{\frac{|s-s'|^{2H}+|t-t'|^{2H}}{2T^{2H}}}-1\right|=0,
 \EQNY
 which coincide with the local variance and correlation behavior of $Z_u(s,t)$ in proof of Theorem \ref{th1} for case $H\geq \frac{1}{2}$. Similarly as in proof of Theorem \ref{th1}, we derive that for $H>\frac{1}{2}$
 \BQNY
  \pk{\sup_{(s,t)\in E_{u,2}} Z_{u,1}(s,t)>m_1(u)}\sim \Psi\left(m_1(u)\right), \quad u\rw\IF;
 \EQNY
 and for $H=\frac{1}{2}$
 \BQNY
 \pk{\sup_{(s,t)\in E_{u,2}} Z_{u,1}(s,t)>m_1(u)}\sim \left(\mathcal{P}_{1/2}^{1}\right)^2\Psi\left(m_1(u)\right),\quad u\rw\IF.
 \EQNY
 Inserting  the above asymptotics and (\ref{upper1}), (\ref{Piterbarg})  in (\ref{main1}), we establish the claim. \\
 {\it \underline{Case $0<H< \frac{1}{2}$}}. Observe that
 \BQNY
 \pk{\sup_{0\leq t\leq T}U_t>u}&=&\pk{\sup_{(s,t)\in A} (X_t-X_s)>u}\\
&=&\pk{\sup_{(s,t)\in A}(B_H(t)-B_H(s)-\frac{1}{2}(t^{2H}-s^{2H})+\mu(t-s))>u}\\
&=&\pk{\sup_{(s,t)\in A} Z_{u,2}(s,t)>m_2(u)},
\EQNY
where $$Z_{u,2}(s,t)=\frac{B_H(t)-B_H(s)}{u-\mu(t-s)+\frac{1}{2}(t^{2H}-s^{2H})}m_2(u), \quad m_2(u)=\inf_{0\leq s\leq T}\frac{u-\mu(T-s)+\frac{1}{2}(T^{2H}-s^{2H})}{(T-s)^H}.$$
Thus we have
\BQN\label{main1}
\pk{\sup_{(s,t)\in E_{u,3}} Z_{u,2}(s,t)>m_2(u)}&\leq& \pk{\sup_{0\leq t\leq T}U_t>u}\nonumber\\
&\leq& \pk{\sup_{(s,t)\in E_{u,3}} Z_{u,2}(s,t)>m_2(u)}+\pk{\sup_{(s,t)\in A\setminus E_{u,3}} Z_{u,2}(s,t)>m_2(u)},
\EQN
where $$E_{u,3}=[0, s_u+(\ln m_2(u))/m_2(u)]\times[T-(\ln m_2(u))^2/(m_2(u))^2,T].$$
In light of Lemma \ref{lem2}, it follows that for $u$ sufficiently large
\BQNY
\sup_{(s,t)\in A\setminus E_{u,3}}\sqrt{Var(Z_{u,2}(s,t))}\leq 1-Q_3\left(\frac{\ln m_2(u)}{m_2(u)}\right)^2,
\EQNY
and direct calculation shows that
\BQNY
\mathbb{E}\left((Z_{u,2}(s,t)-Z_{u,2}(s',t'))^2\right)\leq Q_4(|t-t'|^{2H}+|s-s'|^{2H}), \quad (s,t), (s',t')\in A.
\EQNY
By  Piterbarg Theorem, we have for $u$ sufficiently large
\BQN\label{neg}
\pk{\sup_{(s,t)\in A\setminus E_{u,3}} Z_{u,2}(s,t)>m_2(u)}\leq Q_5 (m_2(u))^{\frac{2}{H}}\Psi\left(\frac{m_2(u)}{1-Q\left(\frac{\ln m_2(u)}{m_2(u)}\right)^2}\right).
\EQN
Next we consider $\pk{\sup_{(s,t)\in E_{u,3}} Z_{u,2}(s,t)>m_2(u)}$. Rewrite
$$\pk{\sup_{(s,t)\in E_{u,3}} Z_{u,2}(s,t)>m_2(u)}=\pk{\sup_{(s,t)\in E_{u,4}} Z_{u,2}(s_u+\Delta_1(u)s,T-\Delta_1(u)t)>m_2(u)}$$
where
$$ E_{u,4}=[-s_u/\Delta_1(u), (\ln m_2(u))/(m_2(u)\Delta_1(u))]\times [0, (\ln m_2(u))^2/((m_2(u))^2\Delta_1(u))], \quad \Delta_1(u)=2^{\frac{1}{2H}}T(m_2(u))^{-\frac{1}{H}},$$
and $s_u$ is defined in Lemma \ref{lem2}.
Lemmas \ref{lem2} and \ref{lem3} lead to
\BQNY
\lim_{u\rw\IF}\sup_{(s,t)\in E_{u,4}}\left|\frac{1-\sqrt{Var(Z_{u,2}(s_u+\Delta_1(u)s,T-\Delta_1(u)t))}}{\frac{H(1-H)}{2T^2}(\Delta_1(u))^2s^2+\frac{H}{T}\Delta_1(u)t}-1\right|=0,
\EQNY
and
\BQNY
\lim_{u\rw\IF}\sup_{(s,t), (s',t')\in E_{u,4}}\left|(m_2(u))^2\frac{1-Corr\left(Z_{u,2}(s_u+\Delta_1(u)s,T-\Delta_1(u)t)), Z_{u,2}(s_u+\Delta_1(u)s,T-\Delta_1(u)t))\right)}{|s-s'|^{2H}+|t-t'|^{2H}}-1\right|=0.
\EQNY
Next we check the conditions of Lemma \ref{Uniform} in Appendix. Following the same notation as in Lemma \ref{Uniform}, we have that 
$$\nu_1=\lim_{u\rw\IF} (m_2(u))^2\frac{H}{T}\Delta_1(u)=2^{\frac{1}{2H}}H\lim_{u\rw\IF} (m_2(u))^{2-\frac{1}{H}}=0,\quad \nu_2=\lim_{u\rw\IF} (m_2(u))^2\frac{H(1-H)}{2T^2}(\Delta_1(u))^2=0,$$
$$y_{1,2}=\lim_{u\rw\IF}m_2(u)\sqrt{\frac{H(1-H)}{2T^2}}\Delta_1(u)(\ln m_2(u))/(m_2(u)\Delta_1(u))=\IF,$$
$$y_{2,1}=0, \quad y_{2,2}=\lim_{u\rw\IF}(m_2(u))^2\frac{H}{T}\Delta_1(u)(\ln m_2(u))^2/((m_2(u))^2\Delta_1(u))=\IF.$$
Moreover, by Lemma \ref{lem2}, $s_u\sim T^{\frac{1}{1-2H}}u^{-\frac{1}{1-2H}}$, which implies that
$$y_{1,1}=-\lim_{u\rw\IF}m_2(u)\sqrt{\frac{H(1-H)}{2T^2}}\Delta_1(u)s_u/\Delta_1(u)=-Q\lim_{u\rw\IF} u^{1-\frac{1}{1-2H}}=0.$$
Thus by case i) in Lemma \ref{Uniform}, we have that
\BQNY
&&\pk{\sup_{(s,t)\in E_{u,4}} Z_{u,2}(s_u+\Delta_1(u)s,T-\Delta_1(u)t)>m_2(u)}\\
&&\quad \sim \left(\mathcal{H}_H\right)^2\sqrt{\frac{2T^2}{H(1-H)}}\frac{T}{H}\int_0^\IF e^{-t^2}dt \int_0^\IF e^{-s}ds(m_2(u))^{-3}(\Delta_1(u))^{-2}\Psi(m_2(u))\\
&& \quad \sim 2^{-\frac{1}{H}-\frac{1}{2}}T^{3H}\sqrt{\frac{\pi}{H^3(H-1)}}\left(\mathcal{H}_{H}\right)^2u^{\frac{2}{H}-3}\Psi(m_2(u)).
\EQNY
Inserting the above asymptotics and (\ref{neg}) into (\ref{main1}) establishes the claim. This completes the proof. \QED
\section{Appendix}
\subsection{Appendix A}
This subsection is devoted to  the proofs of Lemma \ref{lem1}-\ref{lem2}.\\
\prooflem{lem1} Note that for any $\delta>0$ and $u$ sufficiently large, the maximum of $\sigma_u^-(s,t)$ over $0\leq s\leq t\leq T$ is only obtained in $[0,\delta]\times[T-\delta, T]$. Next we consider the variance function $\sigma_u^-(s,t)$ over $[0,\delta]\times[T-\delta, T]$. It follows that
\BQNY
1-\frac{\sigma_u^-(s,t)}{\sigma_u^-(0,T)}&=&1-\frac{|t-s|^H}{u+\mu(t-s)-\frac{1}{2}(t^{2H}-s^{2H})}\frac{u+\mu T-\frac{1}{2}T^{2H}}{T^H}\\
&=&1-\frac{\frac{|t-s|^H}{T^H}}{\frac{u+\mu(t-s)-\frac{1}{2}(t^{2H}-s^{2H})}{u+\mu T-\frac{1}{2}T^{2H}}}\\
&=& \left(1-\frac{|t-s|^H}{T^H}\right)(1+o(1))+\left(\frac{u+\mu(t-s)-\frac{1}{2}(t^{2H}-s^{2H})}{u+\mu T-\frac{1}{2}T^{2H}}-1\right)(1+o(1))\\
&=& \frac{H}{T}(T-t+s)(1+o(1))+\frac{-\mu(T-t+s)+\frac{1}{2}(2HT^{2H-1}(T-t)+s^{2H})}{u+\mu T-\frac{1}{2}T^{2H}}(1+o(1))\\
&=& \left(\frac{H}{T}(T-t)+\frac{H}{T}s+\frac{1}{2u}s^{2H}\right)(1+a(\delta,u)), \quad (s,t)\in [0,\delta]\times[T-\delta, T],
\EQNY
as $\delta$ sufficiently small and $u$ sufficiently large, where $\lim_{\delta\rw 0, u\rw \IF}a(\delta,u)=0$.
The fact that $$\frac{H}{T}(T-t)+\frac{H}{T}s+\frac{1}{2u}s^{2H}>0$$ for $(s,t)\in ([0,\delta]\times[T-\delta, T])\setminus \{(0,T)\}$ implies that the maximum point of $\sigma_u^-(s,t)$ over $0\leq s\leq t\leq T$ is  unique and is $(0,T)$. This completes the proof. \QED\\
\prooflem{lem2} For any $\delta>0$ and $u$ sufficiently large, the maximum of $\sigma_u^+(s,t)$ over $0\leq s\leq t\leq T$ is only obtained in $[0,\delta]\times[T-\delta, T]$. Next we focus on $\sigma_u^+(s,t)$ over $[0,\delta]\times[T-\delta,T]$. For $\delta>0$ sufficiently small and $u$ sufficiently large,
\BQNY
1-\frac{\sigma_u^+(s,t)}{\sigma_u^+(0,T)}&=&1-\frac{|t-s|^H}{u-\mu(t-s)+\frac{1}{2}(t^{2H}-s^{2H})}\frac{u-\mu T+\frac{1}{2}T^{2H}}{T^H}\\
&=&1-\frac{\frac{|t-s|^H}{T^H}}{\frac{u-\mu(t-s)+\frac{1}{2}(t^{2H}-s^{2H})}{u-\mu T+\frac{1}{2}T^{2H}}}\\
&=& \left(1-\frac{|t-s|^H}{T^H}\right)(1+o(1))+\left(\frac{u-\mu(t-s)+\frac{1}{2}(t^{2H}-s^{2H})}{u-\mu T+\frac{1}{2}T^{2H}}-1\right)(1+o(1))\\
&=& \frac{H}{T}(T-t+s)(1+o(1))+\frac{\mu(T-t+s)-\frac{1}{2}(2HT^{2H-1}(T-t)+s^{2H})}{u-\mu T+\frac{1}{2}T^{2H}}(1+o(1))\\
&=& \left(\frac{H}{T}(T-t)+\frac{H}{T}s\right)(1+a_1(\delta,u))-\frac{1}{2u}s^{2H}(1+a_2(\delta,u)), \quad (s,t)\in [0,\delta]\times[T-\delta, T],
\EQNY
where $\lim_{\delta\rw 0, u\rw\IF}a_i(\delta,u)=0, i=1,2. $
If $H\geq \frac{1}{2}$, then
\BQNY
1-\frac{\sigma_u^+(s,t)}{\sigma_u^+(0,T)}
= \left(\frac{H}{T}(T-t)+\frac{H}{T}s\right)(1+a_1(\delta,u)), \quad (s,t)\in [0,\delta]\times[T-\delta, T],
\EQNY
which implies that the maximum point of $\sigma_u^+(s,t)$ is obtained at $(0,T)$ and is unique. For $0<H<\frac{1}{2}$,
\BQNY
1-\frac{\sigma_u^+(s,T)}{\sigma_u^+(0,T)}&=&\frac{H}{T}s(1+a_1(\delta,u))-\frac{1}{2u}s^{2H}(1+a_2(\delta,u))\\
&=&\frac{H}{T}s^{2H}\left(s^{1-2H}(1+a_1(\delta,u))-\frac{1}{2u}(1+a_2(\delta,u))\right)<0,
\EQNY
as $s<\left(\frac{(1+a_2(\delta,u))}{2u(1+a_1(\delta,u))}\right)^{\frac{1}{1-2H}}\sim (2u)^{-\frac{1}{1-2H}}.$ 
This implies that the maximum of $\sigma_u^+(s,T)$ over $[0,T]$ is attained over $(0,\delta)$ for 
$\delta>0$ sufficiently small and $u$ sufficiently large. We denote this point by $s_u$. 
Using the fact that $$\frac{\partial\sigma_u^+(s_u,T)}{\partial s}=\frac{-H(T-s_u)^{H-1}\left(u-\mu(T-s_u)+\frac{1}{2}(T^{2H}-s_u^{2H})\right)-(T-s_u)^H(\mu-Hs_u^{2H-1})}
{(u-\mu(T-s_u)+\frac{1}{2}(T^{2H}-s_u^{2H}))^2}=0,$$
we have that
$$s_u=\left(\frac{u}{T}+\frac{1}{2}T^{2H-1}+\frac{\mu(1-H)}{H}+\frac{1}{2T}s_u^{2H}-\frac{\mu(1-H)}{TH}s_u\right)^{\frac{1}{2H-1}}\sim T^{\frac{1}{1-2H}}u^{-\frac{1}{1-2H}}.$$
Next we show that the maximizer of $\sigma_u^+(s,t)$ is $(s_u,T)$ for $0<H<\frac{1}{2}$ and $u$ sufficiently large.
Observe that
\BQNY
1-\frac{\sigma_u^+(s,t)}{\sigma_u^+(s_u,T)}=-\frac{\sigma_u^+(s,T)-\sigma_u^+(s_u,T)}{\sigma_u^+(s_u,T)}
+\frac{\sigma_u^+(s,T)-\sigma_u^+(s,t)}{\sigma_u^+(s_u,T)}.
\EQNY
Direct calculation gives that, as $u\rw\IF$,
\BQNY
\sigma_u^+(s_u,T)&\sim&\frac{T^H}{u},\\
\sigma_u^+(s,T)-\sigma_u^+(s_u,T)&=&\frac{1}{2}\frac{\partial^2\sigma_u^+(s_u,T)}{\partial^2 s }(s-s_u)^2(1+o(1))\sim \frac{H(H-1)T^{H-2}}{2u}(s-s_u)^2,\\
\sigma_u^+(s,T)-\sigma_u^+(s,t)&=&\frac{\partial\sigma_u^+(s,T)}{\partial t}(T-t)(1+o(1))\sim \frac{HT^{H-1}}{u}(T-t), \quad t\rw T.
\EQNY
Thus we have
\BQNY
1-\frac{\sigma_u^+(s,t)}{\sigma_u^+(s_u,T)}= \frac{H(1-H)}{2T^2}(s-s_u)^2(1+o(1))+\frac{H}{T}(T-t)(1+o(1)), \quad u\rw\IF, |s-s_u|, T-t\rw 0.
\EQNY
The above local behavior implies that the maximizer of $\sigma_u^+(s,t)$ is $(s_u,T)$ for $u$ large and  is unique. This completes the proof. \QED

\prooflem{lem3}
Let $\sigma_{H}(s,t):=\sqrt{\Var(B_H(t)-B_H(s))}$. Observe that
$$\sigma_H(s,t)=|t-s|^H,$$
and
\BQNY
&&1-Corr(B_H(t)-B_H(s),B_H(t')-B_H(s'))\\
&&\quad =
\frac{\mathbb{E}\LT\{((B_H(t)-B_H(s))-(B_H(t')-B_H(s')))^2\RT\}
-(\sigma_H(s,t)-\sigma_H(s',t'))^2}{2\sigma_H(s,t)\sigma_H(s',t')}\\
&&\quad =
\frac{\mathbb{E}\LT\{((B_H(t)-B_H(t'))-(B_H(s)-B_H(s')))^2\RT\}
-(\abs{t-s}^{H}-\abs{t'-s'}^{H})^2}{2\abs{t-s}^{H}\abs{t'-s'}^{H}}\\
&&\quad =
\frac{\abs{t-t'}^{2H}+\abs{s-s'}^{2H}+(|t-s|^{2H}+|t'-s'|^{2H}-|t-s'|^{2H}-|t'-s|^{2H})
-(\abs{t-s}^{H}-\abs{t'-s'}^{H})^2}{2\abs{t-s}^{H}\abs{t'-s'}^{H}}.
\EQNY
Using Taylor formula, we have that for $(s,t)\in [0,\delta_u]\times [T-\delta_u,T]$,
with  $\lim_{u\rw\IF}\delta_u=0$ and $u$ sufficiently large
\BQNY
|t-s|^{2H}-|t-s'|^{2H}-(|t'-s|^{2H}-|t'-s'|^{2H})&=&2H(|\theta_1-s|^{2H-1}-|\theta_1-s'|^{2H-1})(t-t')\\
&=&2H(2H-1)(\theta_1-\theta_2)^{2H-2}(s-s')(t-t'),\\
(\abs{t-s}^{H}-\abs{t'-s'}^{H})^2&=&(H\theta_3(t-t'-s+s'))^2,
\EQNY
where $\theta_1\in (t,t')$, $\theta_2\in (s,s')$ and $\theta_3\in (t-s, t'-s')$. Moreover,
\BQNY
\lim_{u\rw\IF}\lim_{s,t\in [0,\delta_u]\times[T-\delta_u,T]}\left||t-s|^{H}- T^H\right|=0.
\EQNY
Consequently, for $\lim_{u\rw\IF}\delta_u=0$
$$\lim_{u\rw\IF}\sup_{(s,t), (s',t')\in [0,\delta_u]\times [T-\delta_u, T]}\left|\frac{1-Corr\left(B_H(t)-B_H(s), B_H(t')-B_H(s')\right)}{\frac{|s-s'|^{2H}+|t-t'|^{2H}}{2T^{2H}}}-1\right|=0.$$

\subsection{Appendix B}
In this subsection we present some useful results derived in \cite{Uniform2017}.
 First we give an accommodated to our needs version of Theorem 3.2 in \cite{Uniform2017}. Let  $X_u(s,t), (s,t)\in \prod_{i=1,2}[a_i(u), b_i(u)]$ with $0\in \prod_{i=1,2}[a_i(u), b_i(u)]$, be a family of centered continuous
Gaussian random  fields with variance function $\sigma_u(s,t)$ satisfying, {as $u\rw\IF$},
\BQN\label{vvar3}
\sigma_u(0,0)=1, \ \ \sup_{(s,t)\neq (0,0), (s,t)\in \prod_{i=1,2}[a_i(u), b_i(u)]}\left|\frac{1-\sigma_u(s,t)}
{ \frac{|s|^{\beta_1}}{g_1(u)}+\frac{|t|^{\beta_2}}{g_2(u)}}-1\right|
\EQN
with   $\beta_i>0, i=1,2$,  $\lim_{u\rw\IF}g_i(u)=\IF, i=1,2$, $\lim_{u\rw\IF}\frac{|a_i(u)|^{\beta_1}}{g_1(u)}+\frac{+|b_i(u)|^{\beta_2}}{g_2(u)}=0, i=1,2,$
and correlation function satisfying
\BQN\label{ccor3}
\lim_{u\rw\IF}\sup_{(s,t), (s',t')\in \prod_{i=1,2}[a_i(u), b_i(u)], (s,t)\neq (s',t')}\left|n^2(u)\frac{1-Corr(X_u(s,t), X_u(s',t'))}{|s-s'|^{\alpha}+|t-t'|^{\alpha}}-1\right|=0,
\EQN
with $\alpha\in (0,2]$ and
$\lim_{u\rw\IF}n(u)=\IF$.

{We suppose that} $\lim_{u\rw\IF}\frac{n^2(u)}{g_i(u)}=\nu_i\in [0,\IF], i=1,2$.
\BEL\label{Uniform}
Let $X_u(s,t), (s,t)\in \prod_{i=1,2}[a_i(u), b_i(u)]$ with $0\in \prod_{i=1,2}[a_i(u), b_i(u)]$  be a family of centered continuous Gaussian random fields satisfying (\ref{vvar3}) and (\ref{ccor3}).\\
i) If  $\nu_i=0, i=1,2$ and for $i=1,2$, $$\lim_{u\rw\IF}\frac{(n(u))^{2/\beta_i}a_i(u)}{(g_i(u))^{1/\beta_i}}=y_{i,1}, \quad \lim_{u\rw\IF}\frac{(n(u))^{2/\beta_i}b_i(u)}{(g_i(u))^{1/\beta_i}}=y_{i,2},
\quad  \lim_{u\rw\IF}\frac{(n(u))^{2/\beta_i}(a_i^2(u)+b_i^2(u))}{(g_i(u))^{2/\beta_i}}=0,$$
 with $-\IF\leq y_{i,1}<y_{i,2}\leq \IF$,   then
\BQNY
\pk{\sup_{(s,t)\in \prod_{i=1,2}[a_i(u), b_i(u)]}X_u(s,t)>n(u)}
\sim\left(\mathcal{H}_{\alpha/2}\right)^2\prod_{i=1}^{2}\int_{y_{i,1}}^{y_{i,2}}e^{-|s|^{\beta_i}}ds \prod_{i=1}^{2}\left(\frac{g_i(u)}{n^2(u)}\right)^{1/\beta_i}\Psi(n(u)).
\EQNY
ii) If $\nu_i\in (0,\IF)$ and further $\lim_{u\rw\IF}a_i(u)=a_i\in[-\IF,0], \lim_{u\rw\IF}b_i(u)=b_i\in [0,\IF] $, then
 \BQNY
\pk{\sup_{(s,t)\in \prod_{i=1,2}[a_i(u), b_i(u)]}X_u(s,t)>n(u)}
\sim \prod_{i=1}^{2}
\mathcal{P}_{\alpha/2}^{v_{i},\beta_i}([a_i,b_i])\Psi(n(u)),
\EQNY
where $$\mathcal{P}_{\alpha/2}^{v_{i},\beta_i} ([a_i,b_i])=\mathbb{E}\LT\{\sup_{t\in [a_i,b_i]}e^{\sqrt{2}B_{\alpha/2}(t)-|t|^{\alpha}-\nu_i|t|^{\beta_i}}\RT\}\in(0,\IF),\ i=1,2.$$
iii) If $\nu_i=\IF, i=1,2$, then
\BQNY
\pk{\sup_{(s,t)\in\prod_{i=1,2}[a_i(u), b_i(u)]}X_u(s,t)>n(u)}
 \sim \Psi(n(u)).
\EQNY
\EEL
Next we give a simpler version of Proposition 2.2 in \cite{Uniform2017}. Denote by  $\Lambda(u)$ a series of index sets depending on $u$ and by $[a_1,a_2]\times[b_1,b_2]$ a rectangle with $a_1<a_2$ and $b_1<b_2$.  Let  $X_{u,k,l}(s,t), (s,t)\in [a_1,a_2]\times[b_1,b_2], (k,l)\in \Lambda(u)$ be a family of two-dimensional continuous Gaussian random fields with mean $0$ and variance function $1$. There exists $n_{k,l}(u), (k,l)\in \Lambda(u)$ satisfying \BQN\label{thresholds}\lim_{u\rw\IF}\sup_{(k,l), (k',l')\in \Lambda(u)}\left|\frac{n_{k,l}(u)}{n_{k'l'}(u)}-1\right|=0, \lim_{u\rw\IF}\inf_{(k,l)\in \Lambda(u)}n_{k,l}=\IF,
\EQN
such that the correlation function satisfies
\BQN\label{cor3}
\lim_{u\rw\IF}\sup_{(k,l)\in \Lambda(u)}\sup_{(s,t)\neq (s',t'), (s,t), (s',t')\in [a_1,a_2]\times[b_1,b_2]}\left|\left(n_{k,l}(u)\right)^2\frac{1-Corr\left(X_{u,k,l}(s,t), X_{u,k,l}(s',t')\right)}{|s-s'|^{\alpha_1}+|t-t'|^{\alpha_2}}-1\right|=0,\nonumber\\
\EQN
where $\alpha_i\in (0,2], i=1,2$. \\
Then Proposition 2.2 in \cite{Uniform2017} leads to the following result.
\BEL\label{Uniform1} Let  $X_{u,k,l}(s,t), (s,t)\in E, (k,l)\in \Lambda(u)$ be a family of centered two-dimensional continuous Gaussian random fields with variance function $1$. Assume further that (\ref{thresholds})-(\ref{cor3}) hold. Then
$$\lim_{u\rw\IF}\sup_{(k,l)\in \Lambda(u)}\left| \frac{\pk{\sup_{(s,t)\in [a_1,a_2]\times [b_1,b_2]}X_{u,k,l}(s,t)>n_{k,l}(u)}}{\Psi\left(n_{k,l}(u)\right)}
-\mathcal{H}_{\frac{\alpha_1}{2}}([a_1,a_2])\mathcal{H}_{\frac{\alpha_1}{2}}([b_1,b_2])\right|=0$$
\EEL
Finally, we display a lemma concerning the uniform double maximum, a simpler version of Corollary 3.2 in \cite{Uniform2017}.
Let $E_u$ be a family of non-empty compact subset of $\mathbb{R}^2$ and $A_i\subset [0,S]^2, i=1,2$ be two non-empty compact  subsets of $\mathbb{R}^2$. Denote by $\Lambda_0(u)=\{(k_1,l_1,k_2, l_2): (k_i,l_i)+A_i\subset E_u, i=1,2\}$.
Let $n(u)$ and  $n_{k_i,l_i}(u), (k_i,l_i)+A_i\subset E_u $ be a family of positive functions such that
\BQN\label{threshold1}
\lim_{u\rw\IF}\sup_{(k_i,l_i)+A_i\in E_u}\left|\frac{n_{k_i,l_i}(u)}{n(u)}-1\right|=0, i=1,2, \quad \lim_{u\rw\IF} n(u)=\IF.
\EQN
\BEL\label{Uniform2} Let $X_u(s,t), (s,t)\in E_u$ be a family of centered Gaussian random variance $1$ and correlation function satisfying
\BQNY
\lim_{u\rw\IF}\sup_{(s,t)\neq (s',t'), (s,t),(s',t')\in E_u}\left|(n(u))^2\frac{1-Corr(X_u(s,t), X_u(s',t'))}{|s-s'|^{\alpha_1}+|t-t'|^{\alpha_2}}-1\right|=0
\EQNY
Moreover, there exists $\delta>0$ such that for $u$ large enough
$$Corr(X_u(s,t), X_u(s',t'))>\delta-1, (s,t),(s',t')\in E_u.$$
If further (\ref{threshold1}) is satisfied, then there exits $\mathcal{C}>0, \mathcal{C}_1>0$ such that for all $u$ large
\BQNY
\sup_{(k_1,l_1,k_2,l_2)\in \Lambda_0(u), A_i\subset [0,S]^2, A_i\neq \emptyset, i=1,2 }\frac{\pk{\sup_{(s,t)\in (k_1,l_1)+A_1}X_u(s,t)>n_{k_1,l_1}(u), \sup_{(s,t)\in (k_2,l_2)+A_2}X_u(s,t)>n_{k_2,l_2}(u)}}{e^{-\mathcal{C}_1 (F((k_1,l_1)+A_1, (k_2,l_2)+A_2))^{\frac{1}{2}\min(\alpha_1,\alpha_2)}}S^4\Psi(n_{k_1,l_1,k_2,l_2}(u))}\leq \mathcal{C},
\EQNY
where $$F(A,B)=\inf_{s\in A, t\in B}||s-t||, \quad n_{k_1,l_1,k_2,l_2}(u)=\min(n_{k_1,l_1}(u), n_{k_2,l_2}(u)),$$ and
$\mathcal{C}$ and $\mathcal{C}_1$ are independent of $u$ and $S$.
\EEL
\section*{Acknowledgments}
We sincerely thank Professor Enkelejd Hashorva for his encouragement and support to finish this work. Thanks to  the Swiss National Science Foundation Grant 200021-175752/1.
\COM{Partial support from the Swiss National Science Foundation Projects  200020-159246/1 is kindly acknowledged.
K. D\c{e}bicki also acknowledges  partial support by NCN Grant No .}

 	\bibliographystyle{ieeetr}

\bibliography{Drawdownup}
\end{document}